\documentclass[11pt,fullpage]{article}

\newcommand{\ol}{\setlength{\itemsep}{0pt.}\begin{enumerate}}
\newcommand{\eol}{\end{enumerate}\setlength{\itemsep}{-\parsep}}
\newcommand{\ignore}[1]{}
\begin{document}
\begin{center}
{\bf  Van der Waerden Conjecture for Mixed Discriminants  
 }

\vskip 4pt
{Leonid Gurvits }\\
\vskip 4pt
{\tt gurvits@lanl.gov}\\
\vskip 4pt
Los Alamos National Laboratory, Los Alamos , NM 87545 , USA.
\end{center}




\begin{abstract}
We prove that the mixed discriminant of doubly stochastic $n$-tuples of semidefinite hermitian $n \times n$  matrices is
bounded below by $\frac{n!}{n^{n}}$ and that this 
bound is uniquely attained at the $n$-tuple $( \frac{1}{n} I,...,\frac{1}{n} I )$.
This result settles a conjecture posed by R. Bapat in 1989.
We consider various generalizations and applications
of this result.
\end{abstract} 
 
 \newtheorem{ALGORITHM}{Algorithm}[section]
\newenvironment{algorithm}{\begin{ALGORITHM} \hspace{-.85em} {\bf :} 
}%
                        {\end{ALGORITHM}}
\newtheorem{THEOREM}{Theorem}[section]
\newenvironment{theorem}{\begin{THEOREM} \hspace{-.85em} {\bf :} 
}%
                        {\end{THEOREM}}
\newtheorem{LEMMA}[THEOREM]{Lemma}
\newenvironment{lemma}{\begin{LEMMA} \hspace{-.85em} {\bf :} }%
                      {\end{LEMMA}}
\newtheorem{COROLLARY}[THEOREM]{Corollary}
\newenvironment{corollary}{\begin{COROLLARY} \hspace{-.85em} {\bf 
:} }%
                          {\end{COROLLARY}}
\newtheorem{PROPOSITION}[THEOREM]{Proposition}
\newenvironment{proposition}{\begin{PROPOSITION} \hspace{-.85em} 
{\bf :} }%
                            {\end{PROPOSITION}}
\newtheorem{DEFINITION}[THEOREM]{Definition}
\newenvironment{definition}{\begin{DEFINITION} \hspace{-.85em} {\bf 
:} \rm}%
                            {\end{DEFINITION}}
\newtheorem{EXAMPLE}[THEOREM]{Example}
\newenvironment{example}{\begin{EXAMPLE} \hspace{-.85em} {\bf :} 
\rm}%
                            {\end{EXAMPLE}}
\newtheorem{CONJECTURE}[THEOREM]{Conjecture}
\newenvironment{conjecture}{\begin{CONJECTURE} \hspace{-.85em} 
{\bf :} \rm}%
                            {\end{CONJECTURE}}
\newtheorem{PROBLEM}[THEOREM]{Problem}
\newenvironment{problem}{\begin{PROBLEM} \hspace{-.85em} {\bf :} 
\rm}%
                            {\end{PROBLEM}}
\newtheorem{QUESTION}[THEOREM]{Question}
\newenvironment{question}{\begin{QUESTION} \hspace{-.85em} {\bf :} 
\rm}%
                            {\end{QUESTION}}
\newtheorem{FACT}[THEOREM]{Fact}
\newenvironment{fact}{\begin{FACT} \hspace{-.85em} {\bf :} 
\rm}%
                            {\end{FACT}}

\newtheorem{REMARK}[THEOREM]{Remark}
\newenvironment{remark}{\begin{REMARK} \hspace{-.85em} {\bf :} 
\rm}%
                            {\end{REMARK}}
\newcommand{\alg}{\begin{algorithm}} 
\newcommand{\thm}{\begin{theorem}}
\newcommand{\lem}{\begin{lemma}}
\newcommand{\pro}{\begin{proposition}}
\newcommand{\dfn}{\begin{definition}}
\newcommand{\fac}{\begin{fact}}

\newcommand{\rem}{\begin{remark}}
\newcommand{\xam}{\begin{example}}
\newcommand{\cnj}{\begin{conjecture}}
\newcommand{\prb}{\begin{problem}}
\newcommand{\que}{\begin{question}}
\newcommand{\cor}{\begin{corollary}}
\newcommand{\prf}{\noindent{\bf Proof:} }
\newcommand{\ethm}{\end{theorem}}
\newcommand{\elem}{\end{lemma}}
\newcommand{\epro}{\end{proposition}}
\newcommand{\edfn}{\bbox\end{definition}}
\newcommand{\efac}{\bbox\end{fact}}

\newcommand{\erem}{\bbox\end{remark}}
\newcommand{\exam}{\bbox\end{example}}
\newcommand{\ealg}{\end{algorithm}}
\newcommand{\ecnj}{\bbox\end{conjecture}}
\newcommand{\eprb}{\bbox\end{problem}}
\newcommand{\eque}{\bbox\end{question}}
\newcommand{\ecor}{\end{corollary}}
\newcommand{\eprf}{\bbox}
\newcommand{\beqn}{\begin{equation}}
\newcommand{\eeqn}{\end{equation}}
\newcommand{\wbox}{\mbox{$\sqcap$\llap{$\sqcup$}}}
\newcommand{\bbox}{\vrule height7pt width4pt depth1pt}
\newcommand{\qed}{\bbox}
\def\sup{^}
\def\Tp{Tchebyshef polynomial}
\def\Tps{TchebysDeto be the maximafine $A(n,d)$ l size of a code with distance 
$d$hef polynomials}
\newcommand{\rarrow}{\rightarrow}
\newcommand{\larrow}{\leftarrow}
\newcommand{\grad}{\bigtriangledown}

\overfullrule=0pt
\def\setof#1{\lbrace #1 \rbrace}

\section{Introduction}

An $n \times n$ matrix $A$ is called doubly stochastic if it is nonnegative entry-wise
and its every column and row sum to one. The set of $n \times n$ doubly stochastic
matrices is denoted by $\Omega_{n}$.

Let $S_{n}$ be the symmetric group, i.e. the group
of all permutations of the set $\{1,2, \cdots, n\}$.
Recall that the permanent of a square matrix A is defined by
$$ per(A) = \sum_{\sigma \in S_{n}} \prod^{n}_{i=1} A(i,
\sigma(i)). $$
The famous Van der Waerden Conjecture \cite{M78} states that
$$ min_{A \in \Omega_{n}} D(A) = \frac{n!}{n^{n}}$$ and the minimum is
attained uniquely at the matrix $J_{n}$ in which every entry equals
$\frac{1}{n}$.  The ``modern'' attack at the conjecture began in
fifties (of 20th century), see \cite{M78} for some history, and
culminated with three papers \cite{fr}, \cite{Fa81}, \cite{Eg81}.  In
a very technical paper \cite{fr} S. Friedland got very close to the
desired lower bound by proving that $ min_{A \in \Omega_{n}} D(A) \geq
e^{-n}$. \\ In \cite{Fa81} D.I. Falikman proved the lower bound
$\frac{n!}{n^{n}}$ via ingenious and custom-made arguments.  Finally,
the full conjecture was proved by G.P. Egorychev in \cite{Eg81}.  The
paper \cite{Eg81} capitalized on the simple, but crucial, observation
that the permanent is a particular case of the mixed volume or the
mixed discriminant, which we will define below.  Having in mind this
connection, the main inequality in \cite{Fa81} is just a particular
case of the famous Alexandrov-Fenchel inequalities \cite{Al38}.

Let us consider an $n$-tuple
{ \bf A } = (A$_{1}$, A$_{2}$,... A$_{n}$), where A$_{i} = (A_{i}(k,l):
1 \leq k, l \leq n)$ is a complex $n \times n$ matrix $(1 \leq i
\leq n)$.  Then $\det (\sum t_{i}A_{i})$ is a homogeneous
polynomial of degree n in $t_{1},t_{2}, \cdots, t_{n}$.  The number 
\begin{equation}
D({ \bf A }):=D(A_{1}, A_{2}, \cdots, A_{n}) =
\frac{\partial^{n}}{\partial t_{1} \cdots \partial t_{n}} \det
  (t_{1} A_{1} + \cdots + t_{n}A_{n}) 
\end{equation}
is called the mixed discriminant of $A_{1}, A_{2}, \cdots,
A_{n}$.\\

Mixed discriminants were introduced by A.D. Alexandrov 
as a tool to derive mixed volumes of convex sets (\cite{Al38}, \cite{Sch93}). They are also a 3-dimensional
case of multidimemensional Pascal's determinants \cite{pasc}.
 
There exist many alternative ways to define mixed discriminants.  Let $S_{n}$ be the symmetric group, i.e. the group
of all permutations of the set $\{1,2, \cdots, n\}$.  Then the
following identities hold.
\begin{equation}
D(A_{1}, \cdots A_{n})= \sum_{\sigma, \tau  \in S_{n}}
(-1)^{sgn(\sigma \tau}) \prod^{n}_{i=1} A_{i} (\sigma (i), \tau (i)).
\end{equation}
\begin{equation}
D(A_{1}, \cdots A_{n}) = \sum_{\sigma \in S_{n}} \det(A_{\sigma}), 
\end{equation}
where the $i$th column of $A_{\sigma}$ is the $i$th column of
$A_{\sigma(i)}$. \\
\begin{equation}
  D(A_{1}, \cdots A_{n}) = \sum_{\sigma \in
    S}(-1)^{sgn(\sigma)} per(B_{\sigma}),
\end{equation}
where $B_{\sigma}(k,l) = A_{l}(k, \sigma(k))$.  \\
\beqn
M(A_1,...,A_N) = < (A_1  \otimes...  \otimes A_N) V, V >
\eeqn
where the 
$N^{N}$-dimensional vector $V = V(i_1,i_2,...,i_n): 1 \leq i_{k} \leq N,1 \leq k \leq N $ is defined as follows: 

\noindent
$V(i_1,i_2,...,i_n)$ is equal to $(-1)^{sign(\tau)}$ 
if there exists a permutation $\tau \in S_N $,  and equal
to zero otherwise. 

It follows from the definition of the permanent 
that $per(A) = D(A_{1}, \cdots A_{n})$,
where $A_{j} = Diag(A(i,j): 1 \leq i \leq n), 1 \leq j \leq n$.  \\

In a 1989 paper \cite{BA89}  R.B. Bapat defined the set 
$D_{n}$ of doubly stochastic $n$-tuples.  An $n$-tuple ${ \bf A }=(A_{1},
\cdots, A_{n})$ belongs to $D_{n}$ iff the following properties hold:\\

1. \ \  $A_{i} \succeq 0$, i.e. $A_{i}$ is a positive
semi-definite matrix, $1 \leq i \leq n$.\\

2. \ \ $tr A_{i} = 1$ for $1 \leq i \leq n$. \\

3. \ \ $\sum^{n}_{i=1} A_{i} = I$, where $I$, as usual,
stands for the identity matrix.  \\

One of the problems posed in \cite{BA89} is to determine the minimum
of mixed discriminants of doubly stochastic tuples
$$ min_{A \in D_{n}} D(A) = ?$$
Quite naturally, Bapat conjectured that
$$ min_{A \in D_{n}} D(A) = \frac{n!}{n^{n}}$$
and that it is attained uniquely at ${\bf J}_{n} =:(\frac{1}{n} I,...,
\frac{1}{n} I)$.

In \cite{BA89} this conjecture was formulated for real matrices. We
will prove it in this paper for the complex case, i.e. when matrices
$A_{i}$ above are complex positive semidefinite and, thus,
hermitian. (Recall that a square complex $n \times n$ matrix $A = \{
A(i,j) : 1 \leq i,j \leq n \}$ is called hermitian if $A = A^{*} = \{
\overline{A(j,i)} : 1 \leq i,j \leq n \}$. A square complex $n \times
n$ matrix $A$ is hermitian iff $<A x, x > = <x,A x >$ for all $x \in
C^{n}$.)  One of the main tools we will use below are necessary
conditions for a local minimum under semidefinite constraints. It is
very important,in this optimizational context, that the set of $n
\times n$ hermitian matrices can be viewed as an $n^2$-dimensional
real linear space with (real) inner product $<A,B> = : tr(AB) $.

The rest of the paper will provide a proof of Bapat's conjecture. \\
(The lower bound $\frac{n!}{n^{n}}$ for real symmetric doubly stochastic
$n$-tuples was proved in \cite{GS} .)

\section { Basic Facts about Mixed Discriminants}

\noindent
{\bf Fact 1.}  
\begin{equation}
D(X \alpha_{1}A_{1} Y, \cdots, X \alpha_{i}A_{i}Y, \cdots, X
  \alpha_{n}A_{n}Y) = \det(X) \cdot \det(Y) \cdot \prod^{n}_{i=1}
\alpha_{i} \cdot D(A_{1}, \cdots A_{n}).
\end{equation}

\noindent
{\bf Fact 2.}
\begin{equation}
D(x_{i}y_{i}^{\ast}, \cdots, x_{n}y_{n}^{\ast})=\det
(\sum^{n}_{i=1} x_{i}y_{i}^{\ast}).
\end{equation}

Here $x_{i}y_{i}^{\ast}$ is an $n \times n$ complex matrix of rank one,
$x_{i}$ and $y_{i}$ are $n \times 1$ matrices (column-vectors), $y^{\ast}$ is an
adjoint matrix , i.e $y^{\ast} = \overline{y^{T}}$ .\\

\noindent
{\bf Fact 3.}
$$
D(A_{1},\cdots, A_{i-1}, \alpha A + \beta B, A_{(i+1)}, \cdots
A_{n}) = \\
$$
\begin{equation}
\alpha D(A_{1}, \cdots, A_{i-1}, A, A_{i+1}, \cdots
A_{n})\\ + \beta D(A_{1}, \cdots, A_{i-1}, B, A_{i+1}, \cdots, A_{n}).
\end{equation}

\noindent
{\bf Fact 4.}

\noindent
$D(A_{1}, \cdots, A_{n}) \geq 0$ if $A_{i} \succeq 0$
(positive semidefinite), $1 \leq i \leq n$.

This inequality follows, for instance, from the tensor product representation (5).

\noindent
{\bf Fact 5.}

Suppose that $A_{i} \succeq 0, 1 \leq i \leq n$.  Then
$D(A_{1}, \cdots, A_{n}) > 0$ iff for any $1 \leq i_{1} < i_{2}
\cdots < i_{k} \leq n$\\
the following inequality holds: $Rank (\sum^{k}_{j=1} A_{i_{j}}) \geq k$  \cite{Pa85}. \\ 
This fact is a rather direct corollary of the Rado theorem on the rank of intersection of
a matroid of transversals and a geometric matroid, which is a particular case
of the famous Edmonds' theorem on the rank of inrersection of two matroids \cite{gr:lo:sc}.

\noindent
{\bf Fact 6.} \\

$D(A_{1}, \cdots, A_{n}) >0$ if the $n$-tuple $(A_{1}, 
\cdots A_{n})$ is a doubly stochastic.
This fact follows from Fact 5.

\noindent

{\bf Fact 7.} \\
\beqn
D(A_{1}, \cdots A_{i-1}, X, A_{i+1}, \cdots A_{n}) = tr (X
\cdot Q_{i} ) .
\eeqn

where the matrix $Q_{i} =: (\frac{\partial D}{\partial
  A_{i}})^{T}$ ; it follows from the tensor product representation (5) that
if all the matrices $A_{i}$ are hermitian (i.e. $A_{i}
= A_{i}^{\ast})$ then the mixed discriminant $D(A_1,...,A_n)$ is a real number and  $Q_{i} = Q_{i}^{\ast}$ also
$(1 \leq i \leq n)$. \\
 
All previous facts $1-\cdots-7$ are well known 
(see, e.g., \cite{BA89}).  The next, Fact 8, is quite simple, but
seems to be unknown (at least to the author ) .\\

{\bf Fact 8.(Eulerian matrix identity) }\\
\beqn
 \sum^{n}_{i=1} \langle Q_{i} \omega, A_{i}^{*} \omega \rangle =   \sum^{n}_{i=1} \langle A_{i} Q_{i} \omega,  \omega \rangle =
  D(A_{1}, \cdots , A_{n}) \cdot \langle \omega, \omega \rangle,
\eeqn

where $<\omega, \omega>$ stands for the standard inner 
product in $C^{n}$.\\ 

\prf
Consider the following identity : 
\begin{equation}
D(X A_{1} , \cdots, X A_{i}, \cdots, X A_{n}) = \det(X) \cdot D(A_{1}, \cdots ,A_{n}) , 
\end{equation}
where $X$ is real symmetric nonsingular matrix .
Differentiate its left ant right sides respect to this matrix $X$:
\begin{equation}
 \sum^{n}_{i=1} \hat{Q_{i}}^{T} A_{i}^{T} =  \det(X) X^{-1} \cdot D(A_{1}, \cdots ,A_{n}).
\end{equation}

Here $\hat{Q_{i}}^{T}=\frac{\partial D}{\partial \bullet_{i}}$ evaluated at $(X A_{1} , \cdots, X A_{i}, \cdots, X A_{n})$.\\
Puting $X=I$ we get that\\
\begin{equation}
 \sum^{n}_{i=1} Q_{i}^{T}A_{i}^{T} =  D(A_{1}, \cdots ,A_{n}) \cdot I = \sum^{n}_{i=1} A_{i} Q_{i}
\end{equation}
\eprf

\section{ Basic Facts About Minimizers}

\subsection{ Indecomposability}

Following \cite{GS},\cite{GS1}  we call n-tuple $A=(A_{1}, A_{2},
\cdots A_{n})$ consisting of positive semidefinite hermitian matrices
indecomposable if  $Rank (\sum^{k}_{i=1} A_{i_{j}}) > k$ for all  
$1 \leq i_{1} < i_{2}< \cdots< i_{k} \leq n $, where $1 \leq k<n$.  \\

Also, as in \cite{GS},\cite{GS1} associate with a given $n$-tuple ${\bf
  A} = (A_{1},...,A_{n})$ a set $M({\bf A})$ of $n$ by $n$ matrices
$M({\bf A},W)$, where the indices $W$ run over the orthogonal group $O(n)$,
i.e. $WW^{*} =W^{*}W = I$ ,
and for any orthogonal matrix $W$ with columns $\omega_{1},
...,\omega_{n}$ the matrix $M({\bf A},W)$ has as its $(i,j)$ entry
the inner product $<A_{j} \omega_{i}, \omega_{i}>$.  Matrices
$M({\bf A},W)$ inherit many properties of ${\bf A}$, which means that
sometimes we can simplify things by dealing with matrices and not
$n$-tuples.\\
Let us write down a few of these shared properties in the
following proposition.

\pro
\begin{enumerate}
\item
The $n$-tuple ${\bf A}$ is nonnegative (i.e. consists of positive semidefinite hermitian
matrices ) iff for all $W \in O(n)$,
the matrix $M({\bf A}, W)$ has nonnegative entries.
\item
 The $n$-tuple ${\bf A}$ is indecomposable iff for all $W \in
O(n)$, the matrix $M({\bf A},W)$ is fully indecomposable in the sense of \cite{M78}.  
\item
$\sum^{n}_{i=1} a_{i}=I$ iff for all $W \in O(n)$, the
matrix $M({\bf A},W)$ is row stochastic.
\item
$tr(A_{i})=1$, for $1 \leq i \leq n$, iff for all $W \in
O(n)$, the matrix $M({\bf A},W)$ is column stochastic. 
\item
Therefore ${\bf A}$ is doubly stochastic iff for all $W \in
O(n)$, the matrix $M({\bf A},W)$ is doubly stochastic. \\
\end{enumerate}
\epro

The following fact \cite{GS},\cite{GS1}  states that doubly stochastic $n$-tuples can be
decomposed into indecomposable doubly stochastic tuples.

\noindent
{\bf Fact 9.}

Let ${\bf A}=(A_{i},\cdots,A_{n})$ be a doubly stochastic
n-tuple.  Then there exists a partition $C_1 \cup \cdots \cup C_{k}$ of
$\{1, \cdots, n\}$ such that\\ 

1. \ \ For all $1 \leq s \leq k$
$$ dim(Im(\sum_{i \in C_{s}} A_{i}))=: dim X_{s}= |C_{s}|=:c_{s}.$$

2. \ \  The linear subspaces $X_{s} (1 \leq s \leq k)$ are
pairwise orthogonal and $C^{n} = X_{1} \oplus X_{2} \oplus
\cdots \oplus X_{s}$.  \\

It follows from the definition that, in the notation of
Fact 9, the following identity holds : \\
$$D({\bf A}) = D({\bf A}_{1}) \cdots D({\bf A}_{s}) \cdots D({\bf A}_{k}),$$ \\
where ${\bf A}_{s}$ is a doubly stochastic $c_{s}$-tuple
formed of restrictions of matrices $A_{i}(i \in C_{s})$ on the
subspace $X_{s}= Im (\sum_{i \in C_{s}}A_{i}).$\\

\subsection {Fritz John's Optimality Conditions}

Recall that we are to find the minimum of D({\bf
  A}) on $D_{n}$.  The set of doubly stochastic $n$-tuples is
characterized by the following constraints
$$\sum^{n}_{i=1} A_{i} =I, tr A_{i} = 1, A_{i} \succeq 0.$$

Notice that our constrained optimization problem is defined on the
linear space $H=:H_{n} \oplus \cdots \oplus H_{n}$, where $H_{n}$ is
the linear space of $n \times n$ hermitian matrices. We view the
linear space $H_{n}$ as a $n^2$ dimensional real linear space with the
inner product $<A,B>=tr(AB)$. This inner product extends to tuples via
standard summation.  Though a rather straightforward application of John's
Theorem \cite{John48} (see \cite{GS} ) gives the next result, we
decided to include a proof to make this paper self-contained.

\dfn
Consider a doubly stochastic $n$-tuple ${\bf A}=(A_{1}, A_{2}, \cdots, A_{n})$.
Present positive semidefinite matrices $A_i \succeq 0$ in the following block form with respect
to the orthogonal decomposition $C^n = Im(A_i) \oplus Ker(A_i)$ :
\beqn
A_i = \left( \begin{array}{cc}
		  \tilde{A_{i}} & 0\\
		    0 & 0\end{array} \right), A_i \succ 0 ; 1 \leq i \leq n.
\eeqn
Define
a cone of admissible directions as follows:
$$
K_{0}=\{ (Z_{1}, Z_{2}, \cdots, Z_{n}): \mbox{ there exists } \epsilon > 0 \mbox{ such that the tuple } 
$$
$$ 
 (A_{1}+ \epsilon Z_{1},A_{2}+ \epsilon Z_{2}, \cdots, A_{n}+ \epsilon Z_{n})  \mbox { is doubly stochastic.}\}
$$
I.e. $K_{0}$ is a minimal convex cone in the linear space of hermitian $n$-tuples $H=:H_{n} \oplus \cdots \oplus
H_{n}$,  which contains all $n$-tuples $\{{\bf B} -{\bf A} : {\bf B} \in D_{n} \}$.  
We also define the following two convex cones $K_1$,$K_2$  and one linear subspace $K_3$ of $H=:H_{n} \oplus \cdots \oplus H_{n}$ : \\
 
$ K_{1} = \{ (B_{1}, B_{2}, \cdots, B_{n})$, where the matrices $B_i$ are hermitian and
$$
B_i = \left( \begin{array}{cc}
		  B_{i;1,1} & B_{i;1,2}\\
		   B_{i;2,1}   & B_{i;2,2}\end{array} \right) ; Im(B_{i;2,1}) \subset Im(B_{i;2,2}),B_{i;2,2} \succeq 0,  1 \leq i \leq n.
$$

$K_{2} = \{ (B_{1}, B_{2}, \cdots, B_{n})$, where the matrices $B_i$ are hermitian and
$$
B_i = \left( \begin{array}{cc}
		  B_{i;1,1} & B_{i;1,2}\\
		   B_{i;2,1}   & B_{i;2,2}\end{array} \right) ; B_{i;2,2} \succeq 0,  1 \leq i \leq n.
$$

$$
K_{3} = \{ (C_{1}, C_{2}, \cdots, C_{n}): C_i \in H_n, \  tr(C_{i})=0 ( 1 \leq i \leq n ) \mbox{ and } C_{1}+...+C_{n}=0 \}.
$$
\edfn

\pro
\begin{enumerate}
\item
$K_{0} = K_1 \cap K_3 $.
\item
The closure $\overline{K_1} = K_2$.
\item
The closure $\overline{K_0} = K_2 \cap K_3 $.
\end{enumerate}
\epro

\prf 
\begin{enumerate}
\item
Recall that a hermitian block matrix with strictly positive definite block $D_{1,1} \succ 0$
$$
D = \left( \begin{array}{cc}
		  D_{1,1} & D_{1,2}\\
		   D_{2,1}   & D_{2,2}\end{array} \right) 
$$
is positive semidefinite iff $D_{2,2} \succeq 0$ and 
$D_{2,2} \succeq D_{2,1}D_{1,1}^{-1}D_{2,1}^{*}$. This proves that an $n$-tuple $(B_{1}, B_{2}, \cdots, B_{n}) \in K_1$
iff there exists $\epsilon > 0$ such that $A_{i} + \epsilon B_{i} \succeq 0, 1 \leq i \leq n $. Intersection with $K_3$
just enforces the linear constraints $tr(A_i) = 1,  1 \leq i \leq n $ and $\sum_{1 \leq i \leq n} A_i = I$.
\item
This item is obvious.
\item
Clearly, the closure $\overline{K_1 \cap K_3} \subset \bar{K_1} \cap \bar{K_3} = K_2 \cap K_3 $. We need to prove the reverse
inclusion $\bar{K_1} \cap \bar{K_3} = K_2 \cap K_3  \subset \overline{K_1 \cap K_3}$. Consider the following hermitian $n$-tuple
$$
\Delta_i = \frac{1}{n} I - A_i, \Delta_i = \left( \begin{array}{cc}
		  \frac{1}{n} I - \tilde{A_{i}} & 0\\
		    0 & \frac{1}{n} I\end{array} \right) ; 1 \leq i \leq n.
$$
Clearly, the tuple $(\Delta_1,...,\Delta_n)$ is admissible. The important thing is that the $(2,2)$ blocks of matrices $\Delta_i$
are strictly positive definite. Therefore if $(B_1,...,B_n) \in K_2 \cap K_3$ then for all $\epsilon > 0$ the
tuple $(B_1 +\epsilon \Delta_1,...,B_n +\epsilon \Delta_n) \in K_{0} = K_1 \cap K_3$. This proves that 
$\bar{K_1} \cap \bar{K_3} = K_2 \cap K_3  \subset \overline{K_1 \cap K_3}$ and thus that $\bar{K_0} = K_2 \cap K_3 $.
\end{enumerate}
\eprf

\thm If a doubly stochastic $n$-tuple ${\bf
  A}=(A_{1}, A_{2}, \cdots, A_{n})$ is a (local) minimizer then
there exists a hermitian matrix R and scalars $\mu_{i} (1 \leq i
\leq n)$ such that
\beqn
(\frac{\partial D}{\partial A_{i}})^{T} =: Q_{i} = R+ \mu_{i}
I+P_{i}, 
\eeqn
where the matrices $P_{i}$ are positive semi-definite and
$A_{i}P_{i}= P_{i}A_{i} = 0, 1 \leq i \leq n $. 
\ethm

\prf

If a doubly stochastic $n$-tuple ${\bf A}$ is a (local) minimizer then $<Q_{1},Z_{1}>+ \cdots + <Q_{n},Z_{n}> \geq 0$
for all admissible tuples  $(Z_{1}, Z_{2}, \cdots, Z_{n}) \in K_{0} $.
In other words the hermitian $n$-tuple $(Q_{1}, Q_{2}, \cdots, Q_{n})$ belongs to a dual cone $K_{0}^{\prime}$.
It is well known and obvious that a dual cone $\bar{K_{0}}^{\prime}$ of a clossure is equal to $K_{0}^{\prime}$.
By Proposition 3.3 the closed cone $\bar{K_{0}} = K_2 \cap K_3 $, and the convex cones $K_2, K_3$ are closed.
Therefore, see, for instance, \cite{rok}, the closed convex dual cone $\bar{K_{0}}^{\prime} =  K_2 ^{\prime} + K_3^{\prime}$.
We get by a straigthforward inspection that
$$
K_2 ^{\prime} = \{(P_{1}, P_{2}, \cdots, P_{n}) : P_i = \left(
		  \begin{array}{cc} 0 & 0\\ 0 & \tilde{P_i}\end{array}
		  \right), \tilde{P_i} \succeq 0, 1 \leq i \leq n,
$$
$K_3 ^{\prime} = \{(R + \mu_{1}, R + \mu_{2}, \cdots, R + \mu_{n} ) $,
where the matrix $R$ is hermitian and $\mu_{i},1 \leq i \leq n $ are
real.\\ Therefore, we get that $ Q_{i}^{T} = R+ \mu_{i}I+P_{i} $, for
some hermitian matrix $R$, real $\mu_{i}$ and positive semidefinite
$P_{i} \succeq 0$ satisfying the equality $A_{i}P_{i}= P_{i}A_{i} = 0,
1 \leq i \leq n $.
\eprf

\cor
In notations of Theorem 3.4 the following identities
hold:
$$ D({\bf A}) = tr(A_{i} \cdot Q_{i}) = tr(A_{i}(R+\mu I);$$

$$<A_{i} \omega, (R+\mu_{i}I)\omega > =<A_{i} \omega, Q_{i} \omega>,
\omega \in C^{n} ;$$ 

$$ D({\bf A}) = \sum_{1 \leq i \leq n} <A_{i} \omega, (R+\mu_{i}I)\omega > , <\omega ,\omega> = 1 . $$ 
\ecor 
\prf It follows directly from the identity
$A_{i}P_{i}=0$, Facts 7,8 and the hermiticity of all the matrices
involved here.
\eprf

The following simple Lemma will be used in the proof of
uniqueness.
\lem
Let us consider a doubly stochastic $n$-tuple 
$$ {\bf A} = (A_{1}, A_{2}, \frac{1}{n} I \cdots,  \frac{1}{n} I),$$
where $A_{1}, A_{2} \geq 0, tr A_{1} = tr A_{2} = 1$ and $A_{1} +
A_{2} = \frac{2}{n} I$.  Then $D({\bf A}) = \frac {n!}{n^{n}} +
tr ((A_{1} - \frac{1}{n} I) \cdot (A_{1} - \frac{1}{n}
I)^{\ast}))\frac{(n-2)!}{n^{n-2}} $.\\
\elem

\prf First, notice that the matrices in this tuple commute.  Thus,
 $D({\bf A})=per(C_{ij}(1 \leq i,j \leq n)$, where the first column of
 matrix $C$ is equal to $\frac{1}{n} e +\tau$, second column to
 $\frac{1}{n} e -\tau$; all other columns are equal to $\frac{1}{n} e
 $.  Here, as usual, $e$ stands for vector of all ones; the vector
 $\tau$ consists of eigenvalues of $A_{1} - \frac{1}{n}I$. Notice
 that $\sum \tau_{i}=0$.  Using the linearity of the permanent in each
 column we get that
$$per(\alpha_{ij})= \frac {n!}{n^{n}} + Per (B),$$ where the first and
second columns of matrix $B$ are equal to $\tau$ and all others to
$\frac{1}{n} e $. 

An easy computation gives that 
$$ Per B = -2 (\sum_{i<j} \tau_{i}\tau_{j}) \frac
{(n-2)!}{n^{n-2}}.$$
But $0 = (\tau_{i} + \cdots +\tau_{n})^{2} = \tau_{1}^{2} + \cdots +
\tau_{n}^{2} + 2 \sum_{i<j} \tau_{i} \tau_{j}.$ 
Thus 
$$
D({\bf A}) = per(C) = \frac{n!}{n^{n}} +
  \frac{(n-2)!}{n^{n-2}} \cdot (\sum_{i=1}^{n} \tau_{i}^{2}) =
    \frac{n!}{n^{n}} + \frac{(n-2)!}{n^{n-2}} tr ((A_{i} -
      \frac{1}{n}I)(A_{i} - \frac{1}{n}I)^{\ast})\;.
$$ 
\eprf

\section { Proof of Bapat's conjecture }

\thm
\begin{enumerate}
\item
$min_{{\bf A} \in D_{n}} D({\bf A}) = \frac {n!}{n^{n}}.$
\item 
The minimum is uniquely attained at $J_{n} =
(\frac{1}{n}I, \frac{1}{n}I, \cdots \frac{1}{n}I).$
\end{enumerate}
\ethm

\prf
\begin{enumerate}
\item
 To make our proof a bit simpler, we will prove the first part of Theorem 4.1 by
induction.  Assume that the theorem is true for $m < n$.  (Case
$n=1$ is obvious).  Therefore, had a minimizing tuple {\bf A}
decomposed into two tuples $B_{1}$ and $B_{2}$, of dimensions
$m_{1}$ and $m_{2}$ respectively, this would imply, using Fact 9,
that $D(A)=D(B_{1})D(B_{2}) \geq \frac{m_{1}!}{m_{1}^{m_{1}}}
  \frac{m_{2}!}{m_{2}^{m_{2_{2}}}} > \frac {n!}{n^{n}}$. 
The last inequality is clearly wrong as $D(A) \leq D (J_{n})=
\frac {n!}{n^{n}}$.  Thus we can assume that any minimizing tuple
is fully indecomposable.  Now apply  to a minimizing tuple
${\bf A}=(A_{1}, \cdots A_{n})$ Theorem 3.4:  there exists a
Hermitian matrix R and scalars $\mu_{1}, \cdots,\mu_{n}$ such
that $P_{i}=Q_{i}-R-\mu_{i}I \geq 0$ and $A_{i}P_{i}=0$.  From
Corollary 3.5 we get that $D(A)=tr(A_{i}(R+\mu_{i}I)$ and
$D(A)=\sum^{n}_{i=1}<A_{i} \omega, (R+
\mu_{i}I) \omega>$ for any normed vector $\omega \in C^{n}$.

Let $W^{\ast} RW=Diag (\theta_{i}, \cdots,\theta_{n})$ for some
  real $\theta_{i}$ and unitary W.  As in Proposition 3.1, define a doubly
  stochastic matrix $B=M(A,W)$.  i.e. $b_{ij} = <A_{i}
    \omega_{j}, \omega_{j}>.$ Here $\omega_{j}$ is a $j$th column
  of W ( or $j$th eigenvector of R).  Writing identities
  $D(A)=tr(A_{i}(R+\mu_{i}I))$ and $\sum^{n}_{i=1} <A_{i}
  \omega_{j}, (R+ \mu_{i}I)\omega_{j}>=D(A) (1 \leq j \leq n)$ in
  terms of the matrix $B$, we obtain the following  systems of linear equations 
$$\mu_{i} + \sum^{n}_{i=1} b_{ij} \theta_{j} = D(A) (1 \leq i \leq n) ;$$ 
$$\theta_{j} + \sum^{n}_{i=1} B_{ij} \mu_{i} = D(A) (1 \leq j \leq n).$$ 
These equations are, of course, encountered also in the matrix
case, where they led to the crucial London's Lemma  \cite{lon} ,\cite{M78}, 
\cite{Fa81},\cite{Eg81}.  Proceeding
in exactly the same way, we easily deduce that $\mu = (\mu_{1},
\cdots,\mu_{n})$ is an eigenvector of $BB^{T}$ with eigenvalue 1;
$\theta$ is an eigenvector of $B^{T}B$ with eigenvalue 1.

It follows from Proposition 3.1 that $B$ is fully indecomposable,
which is equivalent to the fact that $1$ is a simple eigenvalue for
both $BB^{T}$ and $B^{T}B$ (see e.g. \cite{M78} ). 
Therefore  both $\mu$ and $\theta$ are
proportional to the vector $e$ (all ones).  Thus $\mu_{i} =
\alpha$, $\theta_{i}= \beta$ and $\alpha+\beta = D({\bf A})$. \\
Returning to matrices $R+ \mu_{i}I$ we get that $R + \mu_{i} I =
D({\bf A}) \cdot I$. 

Now comes the punch line, 
a generalization of London's Lemma \cite{lon} to mixed
discriminants:

For a minimizing tuple ${\bf A}=(A_{i} \cdots
A_{n})$ the following inequality holds:
\begin{equation}
(\frac{D(A)}{\partial A_{i}})^{T} =: Q_{i} \succeq R+ \mu_{i}I=D({\bf
  A}) \cdot I.
\end{equation}
Indeed, by Theorem 3.4 $Q_{i} = R+ \mu_{i}I +P_i$ and $P_i \succeq 0$. 
After the inequality (16) is established, we are back to
familiar grounds of the Van der Waerden's conjecture proofs \cite{Fa81},
\cite{Eg81}.  Let us introduce the following notations:
$$  
{\bf A}^{i,j}=(A_1,...,A_i,...,A_i,...,A_n), f_{i,j}({\bf A})=\frac{{\bf A}^{i,j} + {\bf A}^{j,i}} {2}.
$$

I.e., we replace jth matrix in the tuple by the ith one to define ${\bf A}^{i,j}$,
replace ith and jth matrices by their arithmetic average to define $f_{i,j}({\bf A})$.

The celebrated Alexandrov-Fenchel inequalities state that 

1. \ \ $D({\bf A})^{2} \geq D({\bf A}^{i,j}) \cdot D({\bf
  A}^{j,i})$. \\
2. \ \ If $A_{i} \succ 0 (1 \leq i \leq n)$ and $D({\bf A})^{2} = D({\bf A}^{i,j})
\cdot D({\bf A}^{i,j})$ \\
then $A_{i} = \tau A_{j}$ for some positive $\tau$. 

 Suppose that ${\bf A}$ is a minimal doubly stochastic $n$-tuple, i.e. that $D({\bf A}) = min_{{\bf B} \in D_{n}} D({\bf B})$.  
Then for $i \neq j$, $f_{i,j}({\bf A})$ is also minimal doubly
stochastic $n$-tuple.  Indeed, double stochasticity is obvious. 

By the linearity of mixed discriminants in each matrix argument we get that 
$D(f_{i,j} ({\bf A})) = \frac{1}{2} D({\bf A}) + \frac{1}{4}
D({\bf A}^{i,j}) + \frac{1}{4} D({\bf A}^{j,i})$.
Also, 
\beqn D({\bf A}^{i,j}) = tr (A_{i} \cdot Q_{j}), D(A^{j,i}) = tr
(A_{j} \cdot Q_{i})\;. \eeqn  (We used here Fact 7). 

As $Q_{k} \succeq D(A)\cdot I (1 \leq k \leq n)$ from inequality (16) and $tr A_{k} 
\equiv 1$ then $D({\bf A}^{i,j}) \geq D(A)$ and $D(A^{j,i}) \geq D(A).$
Thus, we conclude based on Alexandrov-Fenchel's inequalities that
if {\bf A} is a minimal doubly stochastic n-tuple then $D({\bf
  A}^{i,j}) \equiv D({\bf A}) (i \neq j)$ and $D(f_{i,j} (A)) =
D({\bf A})$.

Additionally, if a minimal tuple ${\bf A} = (A_{1}, \cdots,
A_{n})$ consists of positive definite matrices then $A_{i}=A_{j}
(i \neq j)$.  Indeed $A_{i}=\tau A_{j}$ and $tr A_{i} =
tr A_{j}=1$, thus $ \tau= 1$.  As $A_{1} + \cdots + A_{n} = I$, we
conclude that the only ``positive'' minimal doubly stochastic
n-tuple is $J_{n} = (\frac{1}{n}I, \cdots, \frac{1}{n} I)$.  In
any case, as $D({\bf A}^{i,j}) \equiv D({\bf A})(i \neq j)$, then
$D(f_{i,j}({\bf A})) = D({\bf A})$ for minimal tuples {\bf A}. \\

Define the following iteration on n-tuples:
$${\bf A}_{o} = {\bf A},$$
$${\bf A}_{1} = f_{1,2}(A_{o})$$
$$ \cdots $$
$${\bf A}_{n-1} = f_{n-1,n} (A_{n-2})$$
$${\bf A}_{n} = f_{1,2} (A_{n-1})$$
$$ \cdots. $$\\
It is clear that $${\bf A}_{k} \rightarrow \left( \frac{A_{1} + \cdots
  A_{n}}{n}, \cdots, \frac{A_{1} + \cdots A_{n}}{n} \right).$$

As the initial tuple satisfies $A_{1} + \cdots + A_{n} = I$, then
$A_{K} \rightarrow (\frac{1}{n}I, \cdots \frac{1}{n}I)$.  
(Indeed $ Pr_{1} = f_{1,2},...,Pr_{n-1} = f_{n-1,n}$ are orthogonal projectors in the linear finite-dimensional Hilbert
space of hermitian $n$-tuples. By a well known result, 
$\lim_{k \rightarrow \infty} (Pr_{n-1}...Pr_{1})^{k} = Pr$,
where $Pr$ is an orthogonal projector on the linear subspace $L = Im(Pr_{1}) \cap Im(Pr_{2}) \cap...\cap Im(Pr_{n-1})$.
It is obvious that $L = \{(A_1,...,A_n) : A_1 = A_2 =... = A_n \}$. ) 
As the mixed discriminant is a continuous map from tuples to reals, 
we conclude that
$J_{n} = (\frac{1}{n}I, \cdots \frac{1}{n}I)$ is a minimal tuple
and $min_{{\bf A} \in D_{n}} D({\bf A}) = D(J_{n}) =
\frac{n!}{n^{n}}$.\\

\item
Let us now prove the uniqueness.  We only have to prove that any
minimal doubly stochastic n-tuple consist of positive definite
matrices.  Suppose that not.  Then there exists an integer $k
\geq 0$ such that the tuple ${\bf A}_{k}$ has at least one
singular matrix and the next tuple ${\bf A}_{k+1}$ consists of
  positive definite matrices.\\

Assume without loss of generality that ${\bf A}_{k+1} = f_{1,2}
  ({\bf A}_{k})$ and ${\bf A}_{k} = (A_{1}, A_{2}, A_{3}, \cdots,
  A_{n})$.  \\

Then $\frac{A_{1} + A_{2}}{2} = A_{3} = \cdots = A_{n} =
\frac{1}{n} I.$  From Lemma 3.6 we conclude that \\ 
$$
D({\bf A}_{k}) = D(J_{n}) +  \frac{(n-2)!}{n^{n-2}}\cdot tr ((A_{2} -\frac{1}{n}I) \cdot (A_{2} - \frac{1}{n}I)^{\ast}).
$$

As we already know that $J_{n}$ is a minimal tuple and ${\bf
  A}_{k}$ is also a minimal tuple, thus $A_{1}=A_{2} =
\frac{1}{n}I$.  But at least  one of $A_{1}, A_{2}$
suppose to be singular. We got the desired contradiction.\\
\end{enumerate}
\eprf

\cor
 Let us define the set $D_{n,P}$ of hermitian $n$-tuples as follows
$$
 D_{n,P}= \{ {\bf A}=(A_{1}, A_{2}, \cdots, A_{n}) :  A_{i} \mbox{ is positive semi-definite }, tr(A_{i}) \equiv 1 \mbox{ and }
 \sum^{n}_{i=1} A_{i} = P.
$$
Notice that $tr(P)=n$ necessarily.  If the matrix $P$ is sufficiently
close to the identity matrix $I$ then $ min_{A \in D_{n,P}} D(A) =
\frac{n!}{n^{n}}\det(P)$. \ecor

\prf
It follows from the uniqueness part of Theorem 4.1 that if $P$ is sufficiently close to the identity matrix $I$
then there exists a minimal $n$-tuple ${\bf A}=(A_{1}, A_{2}, \cdots, A_{n}) \in D_{n,P}$ consisting of positive definite
matrices. Very similarly to Theorem 3.4 , 
it follows that $Q_{i}=R + \mu_i I (1 \leq i \leq n ), \mu_i \mbox{ is real and } R \mbox{ is hermitian }$.
It is straightforward to prove that under this condition $D(f_{i,j} (A)) =D({\bf A})$.  
 Also, $f_{i,j} (A) \in D_{n,P}$. 
Indeed, 
\beqn D({\bf A}^{i,j}) = tr(A_{i} \cdot Q_{j})=tr(A_{i} \cdot (R + \mu_j I) \mbox{ and } D({\bf A}) = tr(A_{i} \cdot (R + \mu_i I)\;.\eeqn
Thus, $D({\bf A}^{i,j})=tr(A_{i})(\mu_j-\mu_i)+D({\bf A})=D({\bf A})+(\mu_j-\mu_i).$\\
Therefore, 
$$  
D(f_{i,j} ({\bf A})) = \frac{1}{2} D({\bf A}) + \frac{1}{4}D({\bf A}^{i,j}) + \frac{1}{4} D({\bf A}^{j,i})= D({\bf A}).
$$
As in the proof of Theorem 4.1 , the last equality leads to the
minimality of the tuple $( \frac{P}{n}, \cdots, \frac{P}{n})$. 
\eprf

\section{ Motivations and Connections }

The author came across Bapat's conjecture because of the following
theorem \cite{GS},\cite{GS1}.
\thm
Let us consider an indecomposable $n$-tuple $ {\bf A} =(A_1,...,A_n)$ consisting of positive semi-definite
matrices. Then 

1. There exist a unique vector $\alpha$ of positive scalars $\alpha_i; 1 \leq i \leq n $ with product 
equal to $1$ and a positive definite matrix $S$, such that the 
 $n$-tuple $ {\bf B} =(B_1,...,B_n)$, defined by 
$B_i=\alpha_i S A_i S $, is doubly stochastic.

2. The vector $\alpha$ above is the unique minimum of $det (\sum
t_{i}A_{i})$ on the set of positive vectors with product equal to $1$
and $\min_{x_i > 0, \prod_{i=1}^N x_i=1} \det (\sum x_i A_i) =
(\det(S))^{-2}$.  \ethm

Let us define the following important quantity, the capacity of ${\bf A}$:
$$
Cap ({\bf A}) = \inf_{x_i > 0, \prod_{i=1}^N x_i=1} \det (\sum x_i
A_i). $$

Using Theorem 4.1 and Theorem 5.1 we get the following inequality \cite{GS}, \cite{GS1}:

\begin{equation}
        1 \leq \frac {Cap({\bf A})} {D({\bf A})} \leq \frac {n^{n}}{n!}.
\end{equation}

This last inequality played the most important role
in \cite{GS}, \cite{GS1}.  With many other technical details it led to a
deterministic poly-time algorithm to approximate mixed volumes of
ellipsoids within a simply exponential factor.  In the following
subsection we will use it to obtain a rather unusual extension of
the Alexandrov-Fenchel inequality.

\subsection{Generalized Alexandrov-Fenchel inequalities }

Define $L_n$ as a set of all integer vectors $\alpha = (\alpha_1, \ldots,
\alpha_n)$ such
that $\alpha_i \geq 0$ and $\sum_{i=1}^n \alpha_i=n$.\\

For an integer vector
$$
\alpha = (\alpha_1, \alpha_2, \ldots, \alpha_n): \alpha_i \geq 0, \ \sum
\alpha_i = N \
 $$
we define a matrix tuple
$$
{\bf A}^{(\alpha)} = (\underbrace{A_1, \ldots, A_1}_{\alpha_1}, \ldots,
\underbrace{A_k,
\ldots, A_k}_{\alpha_k}, \ldots, \underbrace{A_n, \ldots, A_n}_{\alpha_n}) \,
$$
i.e., matrix $A_i$ has $\alpha_i$ copies in ${\bf A}^{(\alpha)}$.  For a vector
$x=(x_1,\ldots, x_n)$, we define a monomial $x^{(\alpha)} = x_1^{\alpha_1} \ldots
x_n^{\alpha_n}$.

We will use the notations $M^{(\alpha)}$ for the mixed discriminant $D({\bf A}^{(\alpha)})$ and
$Cap^{(\alpha)}$ for the capacity $Cap({\bf A}^{(\alpha)})$. 

\thm
Consider a tuple $A=(A_1, \ldots, A_n)$ of semidefinite hermitian $n \times n$
matrices.  If vectors
$\alpha$, $\alpha^1, \ldots, \alpha^m$ belong to $L_n$ and
$$
\alpha = \sum_{i=1}^n \gamma_i \alpha^i; \quad \gamma_i \geq 0 \,, \quad
\sum \gamma_i =
1$$
then the following hold inequalities hold:
\beqn
\log ( Cap^{(\alpha)} ) \geq \sum \gamma_i \log ( Cap^{(\alpha^i)} ),
\eeqn
\begin{equation}
\log ( M^{(\alpha)} ) \geq \sum \gamma_i \log ( M^{(\alpha^i)} ) -\log (\frac {n^{n}}{n!}).
\end{equation}
\ethm

\prf
Notice that the inequality (21) follows from (20)  via a direct application
of the inequality (19). It remains to prove (20).\\

First, using the arithmetic/geometric means inequality, we get that 
 $Cap^{(\alpha)}\geq d $ iff
$$
\log ( \det (\sum_{i=1}^n A_i \alpha_i e^{x_i}) ) \geq \sum
\alpha_i x_i + \log d
$$
for all real vectors $x= (x_1, \ldots, x_n)$. 

Now, suppose that
$$
\alpha, \alpha^1, \ldots, \alpha^m \in L_n $$
and
$$
\alpha = \sum_{i=1}^m \gamma_i \alpha^i, \quad \gamma_i \geq 0, \quad
\sum_{i=1}^m
\gamma_i = 1.
$$
Then
$$
\log ( \det ( \sum_{i=1}^n A_i \alpha_i^j e^{x_i}) ) \geq
\langle
\alpha^j, x \rangle + \log ( Cap^{(\alpha^j)} ).$$
Multiplying each of the inequalities above by the corresponding $\gamma_i$
and adding
afterwards we get that
$$
\sum_{j=1}^m \gamma_j \log (\det ( \sum_{i=1}^n A_i \alpha_i^j
e^{x_i} ) )
\geq \langle \alpha, X \rangle + \sum_{j=1}^m \gamma_i \log ( Cap^{(\alpha^j)} ).
 $$
As $\log (\det(X))$ is concave for $X \succ 0$, we eventually get the inequality
$$
\log( Cap^{(\alpha)} ) \geq \sum_{j=1}^m \gamma_j \log ( Cap^{(\alpha^j)} ). 
$$
\eprf

The Alexandrov-Fenchel inequalities  can be written as
$$
\log ( M^{(\alpha)} ) \geq   \frac{\log ( M^{(\alpha^{1})} ) +  \log ( M^{(\alpha^{2})} )}{2}.
$$
Here the vector $\alpha=(1,1,...,1); \alpha^{1}=(2,0,1,...,1) ; \alpha^{2}=(0,2,1,...,1)$.
Perhaps the extra $-\log (\frac {n^{n}}{n!})$ in Theorem 5.2 is just an artifact of our proof?
We will show below that an extra factor is needed indeed.

Define $AF(n)$ as the smallest( possibly infimum ) constant one has to subsract from $\sum \gamma_i \log ( M^{(\alpha^i)} )$
in the right
side of (21) in order to get the inequality. Then by Theorem 5.2 , $AF(n) \leq \log(\frac {n^{n}}{n!}) \cong n$.
We will prove below that $AF(n) \geq n \log(\sqrt{2})$ even for tuples consisting of diagonal
matrices. In this diagonal case the mixed discriminant coincides with the permanent.

Consider the following $N \times N$ matrix with nonnegative entries: 
$$
B = \left( \begin{array}{ccccc}
1 & 1\ddots & \ddots & 0 \\[-0.2cm]
0&1 \ddots & \ddots & \ddots \\[-0.2cm]
0 & \ddots & \ddots \\[-0.2cm]
1 & & & \ddots & 1 \end{array} \right) \, \mbox{ i.e. } B = I+ J \,
$$
where $J$ is a cyclic shift.  Assume  without a ``big" loss of
generality that $N =
2k$.  Define
$$
\alpha^1 = (\underbrace{2,2,2,...,2}_k,0,0,...,0), \quad \alpha^2 =
(\underbrace{0,0,0,...,0}_k, 2,2,
\ldots, 2) \; $$
Then,$e =\frac{\alpha^1 + \alpha^2}{2} , e = (1,1,\ldots, 1) $ and
$per(B^{(e)}) = 2 $, 
$per( B^{(\alpha^1)}) = per( B^{(\alpha^2)}) = 2^k =
2^{\frac{N}{2}}.$
(Here the notation $B^{(\alpha)}$ stands for the matrix having $\alpha_{i}$ copies of $i$th column
of the matrix $B$.)
Thus,
$$
\frac{per(B^{(e)})}{\sqrt{per( B^{(\alpha^1)} \cdot per( B^{(\alpha^2)}}}
= \frac{2}{2^{\frac{N}{2}}} \cong \sqrt{2} ^{-N}. $$
\cnj
$$
\lim_{n \rightarrow \infty} \frac{AF(n)}{n} = 1.$$
\ecnj

\section{Further analogs of van der Waerden conjecture }
\subsection{ $4$-dimensional Pascal's determinants }
The following open question has been motivated by the author's
study of the quantum entanglement \cite{stoc}. 

Consider a block matrix

\beqn 
\rho = \left(\begin{array}{cccc}
		  A_{1,1} & A_{1,2} & \dots & A_{1,n}\\
		  A_{2,1} & A_{2,2} & \dots & A_{2,n}\\
		  \dots &\dots & \dots & \dots \\
		  A_{n,1} & A_{n,2} & \dots & A_{n,n}\end{array} \right),
\eeqn
where each block is a $n \times n$ complex matrix. Define

\beqn
QP(\rho) =:  \sum_{\sigma \in S_n} (-1)^{sign(\sigma)}D(A_{1,\sigma(1)},...,A_{n,\sigma(n)}),
\eeqn
where $D(A_1,...,A_n)$ is the mixed discriminant. If instead of the block form (22), to present 
$\rho = \{\rho(i_1,1_2,i_3,i_4) : 1 \leq i_1,1_2,i_3,i_4 \leq n \}$, e.g. as a $4$-dimensional tensor,
then
\begin{eqnarray}
QP(\rho) & = & \frac{1}{N!}\sum_{\tau_{1}, \tau_{2}, \tau_{3}, \tau_{4} \in S_N}(-1)^{sign(\tau_{1}\tau_{2}\tau_{3}\tau_{4})} \nonumber \\
   & &\prod_{i=1}^N \rho(\tau_{1}(i),\tau_{2}(i),\tau_{3}(i), \tau_{4}(i)).  
\end{eqnarray}

In other words, $QP(\rho)$ is, up to $\frac{1}{N!}$ factor, equal to the  $4$-dimensional Pascal's determinant \cite{pasc}. 

Call such a block matrix $\rho$ doubly stochastic if the following conditions hold:
\begin{enumerate}
\item
$\rho \succeq 0$, e.g. the $n^2 \times n^2$ matrix $\rho$ is positive semidefinite.
\item
$\sum_{1 \leq i \leq n} A_{i,i} = I $.
\item
The matrix of traces $\{ tr(A_{i,j} ) : 1 \leq i,j \leq n \} = I$.
\end{enumerate}
A positive semidefinite block matrix $\rho$ is called separable if
$\rho = \sum_{1 \leq i \leq k < \infty} P_i \otimes Q_i$, where the
matrices $P_i \succeq 0, Q_i \succeq 0 : 1 \leq i \leq k$ ; nonseparable
positive semidefinite block matrices are called entangled.  Notice
that in the block-diagonal case $QP(\rho) = D(A_{1,1},...,A_{n,n})$
and our definition of double stochasticity coincides with double
stochasticity of $n$-tuples from \cite{BA89}. Let us denote the closed
convex set of doubly stochastic $n \times n$ block matrices as
$BlD_{n}$, a closed convex set of separable doubly stochastic $n
\times n$ block matrices as $SeD_{n}$.  It was shown in
\cite{stoc} that $\min_{\rho \in BlD_{n}} QP(\rho) = 0$ for $n \geq
3$ and , on the other hand , $\min_{\rho \in SeD_{n}} QP(\rho) > 0$ for $n
\geq 1$.  \cnj $\min_{\rho \in SeD_{n}} QP(\rho) = \frac{n!}{n^n}.$
\ecnj

It is easy to prove this conjecture for $n=2$ , moreover the following equalities hold : 
$$
\min_{\rho \in BlD_{2}} QP(\rho) = \min_{\rho \in SeD_{2}} QP(\rho) = \frac{2!}{2^2} = \frac{1}{2} .
$$

\subsection{Hyperbolic polynomials} 
The following concept  of hyperbolic polynomials originated in the theory of partial differential equations \cite{gar}.
\dfn
A homogeneous polynomial $p(x), x \in R^m$ of degree $n$ in $m$ real varibles is called hyperbolic in the direction $e \in R^m$ 
(or $e$- hyperbolic) if for any $x \in R^m$ the polynomial 
 $p(x - \lambda e)$  in the one variable $\lambda$ has exactly $n$ real roots counting their multiplicities. We assume in this
paper that $p(e) > 0$.

Denote the ordered vector of roots of $p(x - \lambda e)$ as
$\lambda(x) = (\lambda_{1}(x) \geq \lambda_{2}(x)
\geq... \lambda_{n}(x)) $. It is well known that the product of roots
is equal to $p(x)$. Call $x \in R^m$ $e$-positive ($e$-nonnegative) if
$\lambda_{n}(x) > 0$ ($\lambda_{n}(x) \geq 0$).  Define $tr_{e}(x) =
\sum_{1 \leq i \leq n} \lambda_{i}(x)$.

A $k$-tuple of vectors
$(x_1,...x_k)$ is called $e$-positive ($e$-nonnegative) if $x_i, 1
\leq i \leq k$ are $e$-positive ($e$-nonnegative).  Let us fix $n$
real vectors $x_i \in R^m, 1 \leq i \leq n$ and define the following
homogeneous polynomial: \beqn P_{x_1,..,x_n}(\alpha_1,...,\alpha_n) =
p(\sum_{1 \leq i \leq n} \alpha_i x_i) \eeqn Following \cite{khov}, we
define the $p$-mixed value of an $n$-vector tuple ${\bf X} =
(x_1,..,x_n)$ as 
\beqn M_{p}({\bf X}) = : M_{p}(x_1,..,x_n) =
\frac{\partial^n}{\partial \alpha_1...\partial \alpha_n} p(\sum_{1
\leq i \leq n} \alpha_i x_i) 
\eeqn 
\noindent 
Finally, call an $n$-tuple of real
$m$-dimensional vectors $(x_1,...x_n)$ $e$-doubly stochastic if it is
$e$-nonnegative. $tr_{e} x_i = 1 (1 \leq i \leq n)$ and $\sum_{1
\leq i \leq n} x_i = e $.  Denote the closed convex set of $e$-doubly
stochastic $n$-tuples as $HD_{e,n}$.  \edfn

\xam
Consider the following  homogeneous polynomial $p(\alpha_1,...,\alpha_n) = \det(\sum_{1 \leq i \leq n} \alpha_i A_i)$.
If $A_i \succeq 0 : 1 \leq i \leq n $ and $\sum_{1 \leq i \leq n} A_i \succ 0$ then $p(.)$ is hyperbolic
in  the direction $e$, where $e$ is a vector of all ones. If $\sum_{1 \leq i \leq n} A_i = I$ then $e$-double stochasticity
of $n$-tuple ${\bf X}= (e_1,e_2,...,e_n)$ of $n$-dimensional canonical axis vectors is the same as double stochasticity of
$n$-tuple of matrices ${\bf A} = (A_1,...,A_n)$. Moreover $M_{p}({\bf X}) = D({\bf A})$.
\exam

It was proved in a very recent paper \cite{hyp} that $\min_{{\bf X} \in HD_{e,n}} M_{p}({\bf X}) > 0 $.

\cnj
$\min_{{\bf X} \in HD_{e,n}} M_{p}({\bf X}) = p(e) \frac{n!}{n^n} $
\ecnj

\section{Acknowledgements}

It is my great pleasure to thank my friend and colleague Alex
Samorodnitsky. Without our many discussions on the subject this paper
would have been impossible.  I got the first ideas for the proof of Bapat's
conjecture during my Lady Davis professorship at Technion (1998-1999).  \
Many thanks to that great academic institution.


\begin{thebibliography}{99}

\bibitem{gar}
L. Garding, An inequality for hyperbolic polynomials,
Jour. of Math. and Mech., 8(6): 957-965, 1959.

\bibitem{khov} A.G. Khovanskii, Analogues of the Aleksandrov-Fenchel inequalities for
hyperbolic forms, Soviet Math. Dokl. 29(1984), 710-713.

\bibitem{BA89} 
R. Bapat, Mixed discriminants of positive semidefinite matrices, {\sl
Linear Algebra and its Applications} 126, 107-124, 1989. 

\bibitem{GS} 
L. Gurvits and A. Samorodnitsky,  A deterministic polynomial-time algorithm for
approximating mised discriminant and mixed volume,
{\sl Proc. 32 ACM Symp. on Theory of Computing}, ACM, New York, 2000.

\bibitem{GS1}
L. Gurvits and A. Samorodnitsky, A deterministic  algorithm approximating the mixed discriminant and mixed volume,
and a combinatorial corollary, {\sl Discrete Comput. Geom.} 27: 531 -550, 2002.

\bibitem{fr}
S. Friedland, A lower bound for the permanent of a doubly stochastic
matrix, {\sl Annals of Mathematics}, 110(1979), 167-176.

\bibitem{John48} 
F. John, Extremum problems with inequalities as subsidiary conditions,
{\bf Studies and Essays, presented to R. Courant on his 60th
birthday}, Interscience, New York, 1948.\\

\bibitem{stoc}
L. Gurvits, Classical deterministic complexity of Edmonds' problem and Quantum Entanglement, 
{\sl Proc. 35 ACM Symp. on Theory of Computing}, ACM, New York, 2003.

\bibitem{Al38}
 A. Aleksandrov, On the theory of mixed volumes of convex bodies, IV,
Mixed discriminants and mixed volumes (in Russian), {\sl
Mat. Sb. (N.S.)} 3 (1938), 227-251.


\bibitem{Sch93}
 R. Schneider, {\bf Convex bodies: The Brunn-Minkowski Theory},
Encyclopedia of Mathematics and Its Applications, vol. 44, Cambridge
University Press, New York, 1993.

\bibitem{lon}
D. London, Some notes on the vad der Waerden conjecture,
Linear Algebra and Appl. 4 (1971), 155-160.

\bibitem{Fa81} D. I. Falikman, Proof of the van der Waerden's conjecture on the
permanent of a doubly stochastic matrix, {\sl Mat. Zametki}
29, 6: 931-938, 957, 1981, (in Russian).

\bibitem{gr:lo:sc}
M. Gr\"otschel, L. Lovasz and A. Schrijver, {\bf Geometric Algorithms
and Combinatorial Optimization}, Springer-Verlag, Berlin, 1988.

\bibitem{Eg81}
 G.P. Egorychev, The solution of van der Waerden's problem 
for permanents, {\sl Advances in Math.}, 42, 299-305, 1981. 

\bibitem{Pa85}
 A. Panov, On mixed discriminants connected with positive semidefinite
quadratic forms, {\it Soviet Math. Dokl.} 31 (1985).

\bibitem{pasc}
E. Pascal, Die Determinanten, Teubner-Verlag, Leipzig, 1900.

\bibitem{M78}
H.Minc, Permanents, Addison - Wesley, Reading, MA, 1978.

\bibitem{rok}
R. Tyrrell Rockafellar, Convex analysis, Princeton University Press, 1970.

\bibitem{hyp}
L. Gurvits, Combinatorial and algorithmic aspects of hyperbolic polynomials, 2003 ;
available at http://xxx.lanl.gov/abs/math.CO/0404474.

\end{thebibliography}
\end{document}